\documentclass[reqno]{amsart}
\usepackage{amsmath,amsfonts,amssymb, setspace}

\newtheorem{thm}{Theorem}[section]

\newtheorem{prop}[thm]{Proposition}
\newtheorem{defn}[thm]{Definition}

\newcommand{\R}{\mathbb{R}}

\begin{document}

\title{Huber's theorem for hyperbolic orbisurfaces}
\author{Emily B. Dryden}
\address{Department of Mathematics, Bucknell University, Lewisburg, PA 17837}
\email{ed012@bucknell.edu}
\author{Alexander Strohmaier}
\address{Mathematisches Institut, Beringstrasse 1,  D-53115 Bonn, Germany}
\email{strohmai@math.uni-bonn.de}

\begin{abstract}
\noindent 
 We show that for compact orientable hyperbolic orbisurfaces, the Laplace spectrum
 determines the length spectrum as well as the number of singular points
 of a given order. The converse also holds, giving 
 a full generalization of Huber's theorem to the setting of 
 compact orientable hyperbolic orbisurfaces.
\end{abstract}

\maketitle


\vspace{0.4cm}



\section{Introduction and Result}

 We will be interested in compact hyperbolic orbisurfaces, which are a natural
 generalization of compact hyperbolic Riemann surfaces. By ``hyperbolic'' we will mean that the object is endowed with a Riemannian metric of constant curvature -1.  A hyperbolic orbisurface
 can be viewed as a quotient of the hyperbolic plane by a discrete group of
 isometries which is permitted to include elliptic elements. These elliptic elements give
 rise to conical singularities in the quotient surface.  A local neighborhood of a conical singularity looks like the quotient of a disc by the group generated by the rotation through angle $\frac{2 \pi}{n}$ about the disc's center.  We call such a singularity a \emph{cone point of order $n$}, and say that $\mathbb{Z}_n$ is the \emph{isotropy group} associated to the singularity.  For background on orbifolds
 and the eigenvalue spectrum of the Laplace operator in the orbifold context,
 see \cite{emilythesis}
 and the references therein.
 
 Orbifolds which have the same spectrum of the Laplace operator acting on smooth functions are said to be \emph{isospectral}.  
 It is known that in general, there can be at most finitely many isotropy groups (up
 to isomorphism) in a set of isospectral Riemannian orbifolds that share a uniform
 lower bound on Ricci curvature (see \cite{Stanhopebounds}). 
 In fact, N. Shams, E. Stanhope and D. Webb
 have shown in \cite{Shams} that there exist arbitrarily large (but always finite) isospectral sets which
 satisfy this curvature condition, where each element in a given set has points of
 distinct isotropy.  On the other hand, in \cite{DGGWZ} a spectral invariant is exhibited which, 
within the class of all teardrops and footballs, determines the number and order(s) of the cone point(s). 
Teardrops and footballs are orbisurfaces which are topologically the standard sphere $S^2$; the teardrop 
has one cone point, while the football has a cone point at each of the north and south poles, possibly of 
different orders.  These results lead us to ask whether the spectrum determines the orders of the singular 
points for large classes of orbifolds; we answer this question for the class of compact orientable 
hyperbolic orbisurfaces. 
 
 For compact hyperbolic Riemann surfaces, Huber's theorem says that the Laplace spectrum
 determines the length spectrum and vice versa, where the length spectrum is the
 sequence of lengths of all oriented closed geodesics in the surface, arranged in
 ascending order. There has been recent interest in extending Huber's theorem for
 Riemann surfaces to more general settings. 
 In \cite{ElGrMebook}, J. Elstrodt, F. Grunewald and J. Mennicke prove a version of 
 Huber's theorem for discrete cocompact subgroups of $PSL(2, \mathbb{C})$, while 
 L. Parnovskii \cite{Parnov} states an analog for discrete cocompact subgroups of $SO_+(1,n)$.  
 In both cases, the subgroups are permitted to contain elliptic 
 elements, but it is not shown that the spectrum determines these elements explicitly.  
 We prove that the Laplace spectrum of a hyperbolic orbisurface determines its
 length spectrum as well as the number and orders of the singular points. 
 Our definition of the length spectrum in this case is as follows: if $O$ is a hyperbolic
 orbisurface then there exists a discrete subgroup $\Gamma \subset PSL(2,\R)$ such that
 $O$ is isometric to $\Gamma \backslash \mathbb{H}$, where $\mathbb{H}$ denotes the upper half-plane. 
 The unit tangent bundle
 $T_1 O$ can be defined as the quotient $\Gamma \backslash PSL(2,\R)$.
 The geodesic flow on $T_1 O$ is given by the right action of the one parameter
 group
 \begin{gather*}
  a_t= \left( \begin{array}{cc} e^{t/2} & 0 \\ 0 & e^{-t/2} \end{array} \right).
 \end{gather*}
 An oriented periodic geodesic of length $l$ is by definition a curve
 $\gamma: \R \to T_1 O$ such that $\gamma(t+l)=\gamma(t)$
 and such that $\gamma(t_0) a_t=\gamma(t_0+t)$. Curves $\gamma_1$ and $\gamma_2$ 
 are identified if there is a $t_0 \in \R$ such that $\gamma_1(t)=\gamma_2(t+t_0)$.
 The image of such an orbit under the projection
 $T_1 O \to O=\Gamma \backslash PSL(2,\R)/PSO(2)$ is a closed curve in $O$ which is 
 parametrized by arc-length and which is
 locally length minimizing away from cone points. It may pass through cone points, however.
 The ``length spectrum'' is the set of lengths of all 
 periodic geodesics in the orbisurface counting multiplicities.

 \begin{thm}\label{thm:fullhuber}
  Let $O$ be a compact orientable hyperbolic orbisurface. The Laplace
  spectrum of $O$ determines its length spectrum and the number of cone points
  of each possible order.
  Knowledge of the length spectrum and the number of cone points of each order
  determines the Laplace spectrum.
 \end{thm}

 \noindent Shortly after preparing this manuscript, we learned that 
 P. Doyle and J. P. Rossetti had proven this result independently (cf. \cite{DoRo06}).

\

\noindent \textbf{Acknowledgements.}
This work began at the conference ``Recent developments in spectral geometry'' 
in Blossin, Germany, in November 2004, and we would like to thank the organizers 
for their invitation and Andreas Juhl for interesting discussions.

\section{The Proof}

 The second statement can be proved as for Riemann surfaces (see \cite{emilythesis}). 
 Our proof of the first statement is based on the Selberg trace formula for the wave kernel, that is, 
 the distribution $\mathrm{Tr}\cos (t \sqrt{\Delta - \frac{1}{4}})$, where
 $\Delta$ is the Laplace operator on $O$.
 Let $\{\lambda_n^2\}$ be the sequence of eigenvalues of $\Delta$ and denote as
 usual $r_n^2:=\lambda_n^2-\frac{1}{4}$. Then Selberg's trace formula
 for orbisurfaces (see \cite{Hejhalbook1}, \cite{Iwaniec}) reads
 \begin{eqnarray}\label{eqn:ofldSTF}
  \sum_{n=0}^{\infty} h(r_n) & = & \frac{\mu(F)}{4 \pi}
  \int_{-\infty}^{\infty} r h(r) \tanh (\pi r) dr \nonumber \\
  & & \mbox{} + \sum_{\genfrac{}{}{0pt}{}{\{ P \}}{\text{hyperbolic}}}
  \frac{\ln N(P_c)}{N(P)^{1/2} - N(P)^{-1/2}} g[\ln N(P)] \\
  & & \mbox{}
  + \sum_{\genfrac{}{}{0pt}{}{\{ R \}}{\text{elliptic}}} \frac{1}{2 m(R)
  \sin \theta (R) } \int_{-\infty}^{\infty} \frac{e^{-2 \theta (R)
    r}}{1 + e^{-2 \pi r}} h(r) dr \nonumber, 
\end{eqnarray}
where $h$ is any entire function of uniform exponential type and $h(r)=h(-r)$.
The expression $\mu(F)$ denotes the area of a fundamental domain for $\Gamma$.
The sums are over the conjugacy classes of hyperbolic and elliptic elements in $\Gamma$.  
The norm of the hyperbolic conjugacy class $P$ does not depend on the representative chosen, 
and is denoted by $N(P)$.  We let $P_c$ denote the unique 
primitive hyperbolic conjugacy class such that $P=P_c^l$ with $l \in \mathbb{N}$. 
The function $g$ is the Fourier transform of $h$ and thus is a compactly
supported smooth function.
If $R$ is an elliptic conjugacy class, there exists a unique primitive
elliptic conjugacy class $R_c$ such that $R=R_c^l$.
The integer $m(R)$ denotes the order of $R_c$ and 
$\theta(R)=\frac{\pi l}{m(R)}$ where $1 \leq l \leq m(R)-1$.
We may identify the set of primitive elliptic conjugacy classes $R$ in $\Gamma$
with the set of cone points in $O$ of order $m(R)$. 
The set of hyperbolic conjugacy classes $P$ may be identified with the set of closed periodic 
geodesics in $O$ of length $\ln N(P)$. Since in the literature this last statement is usually proved for
groups without elliptic elements, we include the argument here.
Let $\nu(t)=V a_t$ be a periodic geodesic with period $T$, where $V$ is in 
$\Gamma \backslash PSL(2,\R)$. Choose a representative $v \in PSL(2,\R)$, i.e.
$V=[v]$. Since $\nu$ is periodic with period $T$ there
exists a unique element $\gamma \in \Gamma$ such that $v a_T=\gamma v$.
Since $\gamma$ is conjugate to $a_T$ in $PSL(2,\R)$ it is hyperbolic.
It is now easy to see that another representative $v' \in PSL(2,\R)$
gives rise to a conjugate element $\gamma'$ with norm $e^T$. 
Hence, we have a well defined map from periodic geodesics to conjugacy classes of hyperbolic elements.
Conversely, if $\gamma \in \Gamma$ is a hyperbolic element, then there is an element
$x$ in $PSL(2,\R)$ such that $x^{-1} \gamma x=a_{T}$, for some value $T$. Hence, 
$\gamma(t)=[x] a_t$ is a periodic geodesic with period $T$. It is easy to see
that this geodesic does not depend on the choice of $x$ or on the representative
of $\gamma$ in a conjugacy class. This defines a map from conjugacy classes
of hyperbolic elements in $\Gamma$ to periodic geodesics, and this map is clearly the inverse of the
above map.

The Selberg trace formula in the form (\ref{eqn:ofldSTF}) allows us to give a meaning
to the wave trace $\Phi(t)=\mathrm{Tr}\cos (t \sqrt{\Delta - \frac{1}{4}})$ in the
distributional sense. We may define the functional $\Phi$ on $C^\infty_0(\mathbb{R})$
by 
\begin{eqnarray*}
 \Phi(f):=\frac{1}{2}\sum_{n=0}^{\infty} (\hat f(r_n)+\hat f(-r_n)),
\end{eqnarray*}
where $\hat f$ denotes the Fourier transform of $f$; $\hat{f}$ is known to be entire and of
uniform exponential type. This defines a distribution in $\mathcal{D}'(\mathbb{R})$.
Using (\ref{eqn:ofldSTF}) and the above identifications we obtain
\begin{eqnarray}\label{eqn:ofldwave}
 \mathrm{Tr}\cos ( t \sqrt{\Delta - \frac{1}{4}}) & = &
 -\frac{\mu(F)}{8 \pi} \frac{\cosh (t/2)}{\sinh^2(t/2)} +  \sum_{k =1}^{\infty} \sum_{c \in \mathcal{P}}\frac{l_c}{4\sinh(k l_c/2)}(\delta(|t|-k l_c)) \nonumber \\
& &  + \sum_{x \in \mathcal{C}} \Psi_{m(x)}(t).
\end{eqnarray}
Here the first term has to be understood as the distributional derivative
of the distribution $\frac{\mu(F)}{4 \pi} \frac{1}{\sinh(t/2)}$ defined as a principal value (see \cite{GuillZwor}).

The second term is obtained as follows.  $\mathcal{P}$ is the set of \emph{oriented} primitive closed periodic geodesics; 
if we took $\mathcal{P}$ to be the set of unoriented primitive closed periodic geodesics, then we would need to sum 
over $k \in \mathbb{Z} \setminus 0$.  So we changed from a sum over hyperbolic 
conjugacy classes to a sum over all the iterates of the oriented primitive closed geodesics.  
We let $l_c$ denote the length of $c \in \mathcal{P}$, and we use the fact that 
$$
\frac{1}{\sinh \left(\frac{\ln N(P)}{2}\right)} = \frac{2}{N(P)^{1/2}- N(P)^{-1/2}}.
$$
Finally, the delta function appears when we calculate $g[\ln N(P)]$.  In fact, a direct calculation of $g$ yields 
$$
g[\ln N(P)] = \frac{1}{2} [\delta (\ln N(P) - t) + \delta (\ln N(P) + t)],
$$
but it is not hard to see that $\delta(|t| - \ln N(P))$ contains the same information as the sum 
of the two delta functions above.

In the third summand, the set $\mathcal{C}$ is the set of cone points (which correspond to 
primitive elliptic conjugacy classes), 
and $m(x)$ denotes the order of $x \in \mathcal{C}$.
The functions $\Psi_{m}(t)$ are defined as the Fourier transforms
of the functions
\begin{eqnarray*}
 \psi_{m}(r):=\sum_{l=1}^{m-1} \frac{1}{4 m \sin (\pi l/m)}
 \left(\frac{e^{-2\pi r l/m}}{1+e^{-2 \pi r}}+
 \frac{e^{2\pi r l/m}}{1+e^{2 \pi r}}\right).
\end{eqnarray*}

Since $\psi_m(r)$ is exponentially decaying, the functions $\Psi_{m}(t)$
are real analytic and therefore do not contribute singularities to the
wave trace. Knowing the wave trace
we can read off the lengths of all primitive closed geodesics from its
singular support. This can be done as follows. 
The first term in the wave trace formula is completely determined by the singularity at $t=0$
(which reflects the fact that the volume is spectrally determined).
We subtract it from the right-hand side of (\ref{eqn:ofldwave}). Now take the minimal distance, 
$d$, from $0$ to the singular support of the remaining function; $d$ is the length of a primitive closed geodesic, 
which is contributing to the singularity at $t=d$. 
From the corresponding wave trace
invariant we get the number of primitive closed geodesics contributing to this singularity.
Subtract the contribution of those primitive closed geodesics from (\ref{eqn:ofldwave}), and again take 
the minimal distance from $0$ to the singular
support of the remaining function.  Proceed in this way
to find all lengths and multiplicities of primitive closed geodesics.

We now know that the spectrum determines the finite sum 
\begin{equation}\label{eqn:conesum}
 \sum_{x \in \mathcal{C}} \Psi_{m(x)}(t).
\end{equation}
The functions $\psi_{m}(r)$ behave asymptotically 
like $(2m \sin{\frac{\pi}{m}})^{-1} e^{-\frac{2 \pi r}{m}}$ as $r \rightarrow \infty$; 
from this, we conclude that the $\psi_{m}(r)$ are linearly independent for
different values of $m$.  
This in turn implies that the terms contributing to the sum (\ref{eqn:conesum}) by cone
points of different orders are linearly independent. 
 Thus we can read off the orders of the cone points in $O$; the multiplicity
of a summand corresponding to a particular order tells us 
the number of cone points of that order.


\section{Consequences}

Theorem \ref{thm:fullhuber} has implications for the topology of isospectral hyperbolic orbisurfaces.  
To describe these consequences, we first need to define the Euler characteristic
and state a Gauss-Bonnet theorem for orbisurfaces (see \cite{Thnotes}).

\begin{defn}
Let $O$ be an orbisurface with $s$ cone points of orders
$m_1,\ldots,m_s$.  Then we define the (orbifold) Euler characteristic
of $O$ to be 
\begin{displaymath}
\chi(O) = \chi (X_O) - \sum_{j=1}^s (1-\frac{1}{m_j}),
\end{displaymath}
where $\chi(X_O)$ is the Euler characteristic of the underlying topological space of $O$.
\end{defn}

A Riemannian orbisurface is an orbisurface which is endowed with a Riemannian metric.  
The Gauss-Bonnet theorem gives the usual relationship between topology and geometry for these objects: 

\begin{thm}\label{thm:GaussBonnet}
Let $O$ be a Riemannian orbisurface.  Then
\begin{displaymath}
\int_O KdA = 2 \pi \chi (O),
\end{displaymath}
where $K$ is the curvature and $\chi(O)$ is the orbifold Euler
characteristic of $O$.
\end{thm}

Since the volume of a Riemannian orbisurface is spectrally determined via 
Weyl's asymptotic formula (see \cite{Farsi}), we see that for a Riemannian orbisurface with given curvature,
the spectrum determines the orbifold Euler characteristic.  However, since the
orbifold Euler characteristic involves both the genus of the
underlying surface and the orders of the cone points in the
orbisurface, it is not immediately clear that the spectrum determines
the genus.  In the case of compact orientable hyperbolic orbisurfaces, Theorem \ref{thm:fullhuber} says that the spectrum determines the orders of the cone points, and thus the genus.  This observation proves

\begin{prop}
Isospectral compact orientable hyperbolic orbisurfaces have the same underlying topological space.
\end{prop}


\bibliographystyle{amsplain}

\newcommand{\noopsort}[1]{} \newcommand{\printfirst}[2]{#1}
  \newcommand{\singleletter}[1]{#1} \newcommand{\switchargs}[2]{#2#1}
  \def\cprime{$'$}

\end{document}